\documentclass{amsart}

\usepackage{amssymb}
\usepackage[all]{xy}
\usepackage{hyperref}

\usepackage{enumitem}   

\setlist[enumerate]{itemsep=.2em,topsep=.2em,leftmargin=1.25em,itemindent=2.0em}

\newtheorem{thm}{Theorem}

\newtheorem{lem}[thm]{Lemma}
\newtheorem{cor}[thm]{Corollary}

\newtheorem{prop}[thm]{Proposition}

\theoremstyle{definition}
\newtheorem{defn}[thm]{Definition}

\newtheorem{say}[thm]{}
\newtheorem{exmp}[thm]{Example}
\newtheorem{exmps}[thm]{Examples}

\newtheorem{rems}[thm]{Remarks} 

\newtheorem*{ack}{Acknowledgments}      

\newtheorem{defn-thm}[thm]{Definition--Theorem}  
\newtheorem{defn-lem}[thm]{Definition--Lemma}  

\theoremstyle{remark}

\renewcommand{\c}[0]{{\mathbb C}}  

\renewcommand{\o}[0]{{\mathcal O}} 
\newcommand{\z}[0]{{\mathbb Z}}

\renewcommand{\a}[0]{{\mathbb A}}

\newcommand{\p}[0]{{\mathbb P}}

\newcommand{\q}[0]{{\mathbb Q}}
\newcommand{\map}[0]{\dasharrow}
\newcommand{\qtq}[1]{\quad\mbox{#1}\quad}
\newcommand{\spec}[0]{\operatorname{Spec}}
\newcommand{\pic}[0]{\operatorname{Pic}}

\newcommand{\mult}[0]{\operatorname{mult}}

\newcommand{\supp}[0]{\operatorname{Supp}}    
\newcommand{\red}[0]{\operatorname{red}}

\newcommand{\aut}[0]{\operatorname{Aut}}

\newcommand{\ex}[0]{\operatorname{Ex}}

\newcommand{\chr}[0]{\operatorname{char}}

\newcommand{\cl}[0]{\operatorname{Cl}}

\newcommand{\ndeg}[0]{\operatorname{ndeg}}

\newcommand{\simq}[0]{\sim_{\q}}

\newcommand{\tsum}[0]{\textstyle{\sum}}

\def\into{\DOTSB\lhook\joinrel\to}

\def\loccoh#1.#2.#3.#4.{H^{#1}_{#2}(#3,#4)}

\DeclareMathAlphabet{\mathchanc}{OT1}{pzc}%
                                {m}{it}

\newcommand{\gm}[0]{{\mathbb G}_m}

\usepackage[all]{xy}\xyoption{dvips}

\begin{document}
\bibliographystyle{amsalpha}

\hfill\today

 \title[Log K3 surfaces]{Log K3 surfaces with irreducible boundary}
        \author{J\'anos Koll\'ar}

         \begin{abstract}  We determine the automorphism group of an open log   K3 surface with irreducible boundary. 
        \end{abstract}

        \maketitle

   A  {\it log K3} surface over a field $k$ is a pair  $(S, N)$, where
   $S$ is a projective surface over $k$ such that $S_{\bar k}$ is rational, and $N\in |-K_S|$ is a curve with only nodes as singularities.  These are frequently  called {\it anticanonical pairs} or {\it log Calabi-Yau} surfaces.
   We allow  $S_{\bar k}$ to have Du~Val singularities away from $N$.
   The field $k$ can be arbitrary, but some minor issues are  left to the reader in characteristic 2. 
   
   Our aim is to understand whether 
the open surface $S\setminus N $ determines the pair
  $(S, N)$.
 As  a first counterexample, 
let $s\in N$ be a node that is a smooth point of $S$, and $\pi:S'\to S$ its blow up with exceptional curve $E$. Then $(S', N'+E)$ is another log K3 surface, and
$S\setminus N \cong S'\setminus (N'+E)$.
Less obvious examples are given in (\ref{s.b.exmp}.1--4). 
These suggest that the best one can hope for is the following.

\begin{thm} \label{main.thm.1}  Let $k$ be a field and  $(S_i, N_i)$  smooth, log K3 surfaces  over $k$ with geometrically irreducible boundaries.
  Then
  $$
  S_1\setminus N_1 \cong S_2\setminus N_2\quad\Leftrightarrow\quad 
  (S_1, N_1)\cong (S_2, N_2).
  \eqno{(\ref{main.thm.1}.1)}
  $$
\end{thm}

The natural approach to proving (\ref{main.thm.1}.1) would be to show that 
every isomorphism
$\phi^\circ: S_1\setminus N_1 \cong S_2\setminus N_2$
extends to an isomorphism
$\phi: S_1\cong S_2$. This  turns out to hold if $(N_i^2)\leq 2$, but   not  otherwise.

The failure of this isomorphism-extension property is measured by the
difference between  $\aut(S\setminus N)$ and  $\aut(S,N)$.
We show  that there is a natural normal subgroup
$\aut^\circ(S\setminus N)\subset \aut(S\setminus N)$ such that
$$
\aut(S\setminus N)\cong  \aut^\circ(S\setminus N)\rtimes \aut(S,N).
\eqno{(\ref{main.thm.1}.2)}
$$
It is easy to check that  $\aut^\circ(S\setminus N)$ is trivial if either
$N$ is smooth or $(N^2)\leq 0$; see (\ref{must.node.lem}) and (\ref{main.thm.1.pf.1}).
We are thus left with the cases when  $S$ is a  Del~Pezo surface (\ref{a.dp.def}) and
$N$ is a nodal curve. We determine $\aut^\circ(S\setminus N)$ for all  of them.

\begin{thm}  \label{main.thm.2}    Let $k$ be a field,
 $S$  a   Del~Pezzo surface  over $k$ and  $N\in |-K_S|$  a geometrically  irreducible  curve with a node $s\in N$, such that  $S$ is  smooth at  $s$.
    Then $\aut^\circ(S\setminus N)$ is
  \begin{enumerate}
  \item  trivial if $\deg S=1,2$,
  \item $\z/2$ if $\deg S=3$,
  \item trivial if $4\leq \deg S \leq 9$ and the node is unsplit (\ref{node.defn}), and
  \item the  infinite dihedral group  $D_{\infty}$  (\ref{inf.d.defn}) if  $4\leq \deg S \leq 9$ and the node is split.
  \end{enumerate}
\end{thm}

{\it Previous results.}
\cite[4.2.2]{MR0769779} describes an infinite order element of  $\aut\bigl(\p^2\setminus(\mbox{nodal cubic})\bigr)$. The same construction is 
used in \cite[Sec.6]{MR1942244} to obtain highly cuspidal rational curves in $\p^2$.  It turns out  that a  variant 
gives the whole $\aut^\circ(\p^2\setminus N)$; most likely  both authors   suspected this. Versions for other Del~Pezzo surfaces are considered in \cite[Sec.2.1]{mcduff2024singular}.
A more detailed discussion of the relationship with the maps considered in 
\cite{MR0769779, MR1942244, mcduff2024singular} is given in
(\ref{a.giot}).

\begin{rems} If $S$ is a  Del~Pezzo surface, then $\aut(S,N)$ is always finite, and trivial in most cases; see (\ref{nodal.p.say}) for some examples.

  By contrast, $\aut(S,N)$ can be infinite if $(N^2)\leq 0$.
A comprehensive study  is given in
\cite{li2022cone}, which builds on the works of
\cite{MR3314827, MR3415066, friedman2016geometry}.

I do not know how to write down a priori the quotient map
$\aut(S\setminus N)\to \aut(S,N)$.
Thus we determine the structure of $\aut(S\setminus N)$ in all cases, and then observe that they have a   semidirect product structure as claimed in 
(\ref{main.thm.1}.2).
\end{rems}


\begin{ack} The main impetus came from trying to understand the algebraic geometry side of the paper \cite{mcduff2024singular}. 
  I thank D.~McDuff and K.~Siegel for many  e-mails answering my questions and posing new problems.
The many comments and corrections of   P.~Hacking, J.~Li and  K.~Oguiso were also very   helpful.
  
  Partial  financial support    was provided  by  the NSF under grant number
DMS-1901855.
\end{ack}

\section{Outline of the proofs}

For Theorem~\ref{main.thm.2},
the easy cases are (\ref{main.thm.2}.1--3); these are settled in
(\ref{main.thm.1.pf.1}--\ref{main.thm.1.pf.2}).
In order  to prove  (\ref{main.thm.2}.4)
we   first construct some automorphisms.

\begin{say}[Geiser-type involutions]\label{gei.i.i}
   Let
 $S$ be  a   Del~Pezzo surface  of  degree $\geq 4$, and  $N\in |-K_S|$  a geometrically  irreducible  curve with a split node $s\in N$, such that  $S$ is  smooth at  $s$.
   In (\ref{a.obu})   we construct  two  Geiser-type involutions
   $$\sigma_+, \sigma_-\in \aut(S\setminus N),
   \eqno{(\ref{gei.i.i}.1)}
   $$
   indexed by the local branches  $B_{\pm}$ of $N$ at $s$.
There is no natural way to distinguish the branches from each other, so
  it is best to think of $\{\sigma_+, \sigma_-\}$ as a pair of involutions.
  Then we define
  $$
 \aut^\circ(S\setminus N):= \langle \sigma_+, \sigma_-\rangle.
  \eqno{(\ref{gei.i.i}.2)}
   $$
  We check in Paragraph~\ref{bounds.i.say} that
  $\aut^\circ(S\setminus N)\cong D_{\infty}$, and we see in
  (\ref{sd.prod.cor}) that its normalizer contains $\aut(S,N)$.

It  takes longer to show that the  Geiser involutions and $\aut(S,N)$ generate
$\aut(S\setminus N)$.
The approach is modeled on the Noether-Fano method, as implemented by Segre and Manin for cubic surfaces; see \cite{MR0008171, MR0460349} or \cite[Sec.2.1]{ksc}.

The first step is  a computation of how the basic numerical invariants of a cuspidal curve transform under $\sigma_+, \sigma_-$.  The following is equivalent to \cite[2.2.6]{mcduff2024singular}; in the first version  the second case in (\ref{cusp.transf.cor}.2) was  overlooked.
The formulas are simplest using the normalized degree
$$
\ndeg (C):=\tfrac1{\deg S}(-K_S\cdot C)\qtq{and} \mult(C):=\mult_s(C).
\eqno{(\ref{gei.i.i}.3)}
$$
\end{say}

\begin{lem}  \label{cusp.transf.cor}
  Let 
  $S$  be a   Del~Pezzo surface  of degree $d\geq 4$, and  $N\in |-K_S|$  an irreducible  curve with a  split node $s\in N$, such that  $S$ is  smooth at  $s$.
  Let $C\subset S$ be an irreducible curve with only single, cuspidal  point on $N$.  Then, for a suitable choice of the  local branch  $B_{+}$ of $N$,
  $$
  \begin{array}{ccl}
 (\ref{cusp.transf.cor}.1)\qquad   \ndeg\bigl(\sigma_+(C)\bigr)&=&(d{-}1)\ndeg(C)-\mult(C),\qtq{and}\\[1ex]
  (\ref{cusp.transf.cor}.2)\qquad  \ndeg\bigl(\sigma_-(C)\bigr)&=&
    \left\{
      \begin{array}{l}
        \mult(C) -\ndeg(C)\qtq{if $\tfrac{d}{d{-}1}\geq \tfrac{\mult C}{\ndeg C}$, and}\\[1ex]
        (d{-}1)\ndeg(C)-(d{-}2)\mult(C)\qtq{otherwise.}
        \end{array}
      \right.
  \end{array}
  $$
\end{lem}

The proof is given in (\ref{cusp.transf.cor.pf}).

\begin{say}[Bounds on the multiplicity]\label{bounds.i.say}
  Note first that  $2\cdot\mult(C)\leq (N\cdot C)=d\cdot \ndeg(C)$, thus, if
  $d\geq 4$ then 
  $$
  \ndeg\bigl(\sigma_+(C)\bigr)> \ndeg (C),
   \eqno{(\ref{bounds.i.say}.1)}
   $$
   except when $d=4$ and $\mult(C)=2\cdot\ndeg(C)$. 
   Thus we can  increase the degree by one of 
   $\sigma_+,\sigma_-$, hence they   generate an infinite dihedral subgroup of
   $\aut^\circ(S\setminus N)$.

By explicit computation we also see that
$$
\ndeg\bigl(\sigma_-(C)\bigr)<\ndeg(C)
  \eqno{(\ref{bounds.i.say}.2)}
  $$ provided
  $$\ndeg(C)<\mult(C)<2 \cdot\ndeg (C).
  \eqno{(\ref{bounds.i.say}.3)}
  $$
  The upper bound $\mult(C)<2\cdot \ndeg (C)$ holds for almost all curves, see
  Corollary~\ref{q.2m.cor}.
  The exceptional cases are when  $d\in\{4,5\}$ and $\mult(C)=2=2\cdot \ndeg (C)$.

Therefore, using Lemma~\ref{cusp.transf.cor} we can lower the degree of $C$, and get the following.
\end{say}

\begin{cor}  \label{cusp.mult.thm}    Let $(S,N)$ and $C\subset S$ be as in
  Lemma~\ref{cusp.transf.cor}. 
  Then there is a $\phi\in \langle \sigma_+, \sigma_-\rangle$    such that
\begin{enumerate}
\item either $\mult \bigl(\phi(C)\bigr)\leq  \ndeg\bigl(\phi(C)\bigr)$, 
         \item  or $d\in\{4,5\}$, $C$ is rational,  and $\phi(C)\in |-K_S|$ has a $(2,3)$ cusp at $s$.
\end{enumerate}

Moreover, for (\ref{cusp.mult.thm}.1), 
\begin{enumerate}\setcounter{enumi}{2}
\item either $\phi$ is unique and the inequality in (\ref{cusp.mult.thm}.1) is strict,
\item or there are two such  $\phi$  and equality holds in (\ref{cusp.mult.thm}.1) for both.
  \end{enumerate}
\end{cor}

The  inequality  (\ref{cusp.mult.thm}.1) is exactly the one needed for
the  Noether-Fano method.
(The  Noether-Fano method usually works with  the base points of birational transforms of very ample linear systems. In the current situation, studying the cusps of curves works better.)

To complete the proof of  Theorem~\ref{main.thm.2}
we start with an elliptic curve $E\subset S$ that meets $N$ at a single smooth point, see (\ref{E.1pt.lem}).
Given $\psi\in \aut(S\setminus N)$, we apply Corollary~\ref{cusp.mult.thm} to
$\psi(E)$.
Since $E$ is not rational, we are in case (\ref{cusp.mult.thm}.1).
Thus we have a $\phi\in \langle \sigma_+, \sigma_-\rangle$
such that
$$
\mult \bigl(\phi\circ \psi(E)\bigr)\leq  \ndeg\bigl(\phi\circ \psi(E)\bigr).
$$
The Noether-Fano method now says that  $\phi\circ \psi$ is an
automorphism of  $(S, N)$, as needed; see Section~\ref{n-f.meth.sec} for details.\qed

\medskip

Another application  to 
log K3 surfaces
is the study of
 curves $C\subset S\setminus N $
 that are isomorphic to $\a^1$.
 These  are called {\it affine lines} on  $S\setminus N $.
A series of papers 
\cite{Magill_2023,  magill2023staircase, mcduff2024singular} studies those affine lines
$C\subset S$  over $\c$ that are `meaningful' in symplectic geonetry. That is,  a certain Fredholm index is 0; see  \cite[Sec.1.2]{mcduff2024singular}.
These are the affine lines 
for which the birational transform of $C$ on the minimal log resolution of $(S, N+C)$ is a $(-1)$-curve. \cite{Magill_2023,  magill2023staircase, mcduff2024singular} give a complete answer for $\p^2$ and for its 1-point blow-up.

Combining Corollary~\ref{cusp.mult.thm}  with  the multiplicity inequality of \cite{MR1016092} gives  all  affine lines  for $S=\c\p^2$. 

\begin{cor} \label{thm.aff.lines} The affine lines $\a^1\cong C\subset \bigl(\c\p^2\setminus(xyz=x^3+y^3)\bigr)$ form 3 orbits under $\aut\bigl(\c\p^2\setminus(xyz=x^3+y^3)\bigr)$.
  Representative curves in  these orbits are:
  \begin{enumerate}
  \item  $(x=0)$, tangent to a branch of $N$ at the node,
  \item $(3x+3y+z=0)$, flex tangent of $N$, and
  \item  $(21x^2-22xy+21y^2-6xz-6yz+z^2=0)$, conic 6-tangent to $N$.
  \end{enumerate}
\end{cor}

The  proof is given in Section~\ref{al.p2.sec}.
All these affine lines were known to \cite{ MR1942244, mcduff2024singular}, and
the series (\ref{thm.aff.lines}.1)
gives the lines  that have Fredholm index 0.
The proof relies  on  \cite{MR1016092}, which is known only over $\c$, and does not seem to have an analog for other  Del~Pezzo surfaces.

On  the 1-point blow-up of $\p^2$ there are infinitely many orbits of lines with  
Fredholm index 0  \cite{Magill_2023,  magill2023staircase, mcduff2024singular}; and also infinitely many orbits of other lines, see (\ref{a.l.F1.exmp.2}).

\begin{say}[Related works]
Isomorphisms between affine surfaces of the form  $\p^2\setminus(\mbox{curve})$ are considered in
\cite{MR0769779}, giving  a fairly detailed classification.
Remaining cases involving some high degree curves are studied in 
\cite{MR2561193}.

Affine lines in $\a^2=\p^2\setminus(\mbox{line})$ are classified in 
\cite{MR0379502}, and affine lines in
$\p^2\setminus(\mbox{smooth conic})$ are described in \cite{de-du}.
In these cases, a  suitable automorphism transforms any affine line into a line (resp.\ a line or a conic). 

The papers \cite{takahashi1996curves, MR1844627} consider curves in 
$\p^2\setminus(\mbox{smooth cubic})$
whose normalization is $\a^1$; with strongest answers in degrees $\leq 8$. A more general counting problem involving smooth cubics and rational curves 
is considered in    \cite{MR2435425}.

Using birational transformations to study cuspidal curves has been a recurring theme; see for example \cite{MR1016092, MR1457728, MR1942244}.
\cite{MR2280130} gives a topological classification of  rational curves in $\p^2$ that have a unique singularity,  that is a cusp with one Puiseux  pair.
\end{say}

\begin{exmps}\label{s.b.exmp} These examples show that the assumptions   in Theorem~\ref{main.thm.1} are mostly necessary.

  (\ref{s.b.exmp}.1)
  A boundary $N=\cup_iN_i$ could be called minimal if none of its irreducible components are $(-1)$-curves.   If $(N_i^2)\leq -2$ for every $i$ then it is easy to see that
$S\setminus N$ determines $(S, N)$. 
However, 
  Theorem~\ref{main.thm.1} does not extend to all minimal boundaries.
  For example, $\gm\times \gm$ can be compactified to
$\bigl(\p^2, (xyz=0)\bigr)$ and to 
$\bigl(\p^1_x\times \p^1_y, (x_0x_1y_0y_1=0)\bigr)$, neither boundary contains a curve with $-1$ self-intersection.

  (\ref{s.b.exmp}.2)
A more sophisticated version is the following. 
  Let $S$ be a nontrivial Severi-Brauer surface and
  $N\subset S$ the union of 3 conjugate lines.
Note that $N$ is irreducible, but not  geometrically irreducible.

  Blowing up the 3 nodes of $N$ and then contracting the birational transform of $N$ gives the dual Severi-Brauer surface $(S', N')$.
  Thus $S\setminus N\cong S'\setminus N'$, but $S\not\cong S'$.

  (\ref{s.b.exmp}.3) Let   $(S, N)$ be a log K3 surface such that $d:=(N^2)\leq 2$.
  Assume that  $S$ is  smooth at the node  $s\in N$.
  Blowing up $s$, the birational transform $\tilde N$ of $N$ is a smooth rational curve with self-intersection $d-4$. Contracting it we get a
  log K3 surface  $(S', N')$ such that
  $S\setminus N\cong S'\setminus N'$ but $S\not\cong S'$.
  Here $S'$ is singular at the node $s'\in N'$.
  If $d=2$ then $S'$ is also a Del~Pezzo surface with an $A_1$ singularity at 
   the node $s'$.

(\ref{s.b.exmp}.4) If $(N^2)<0$ then $N$ is (analytically) contractible
  $(S,N)\to (S^*, s^*)$. Then $S^*$ is the one-point compactification of
  $S\setminus N$, and $S$ is the minimal resolution of the singularity at $s^*$.
  Thus the isomorphism-extension property holds.
  The same holds if $N$ is reducible with a negative definite intersection matrix.
   
\end{exmps}

\section{Notation and definitions}

 \begin{say}[Del~Pezzo surfaces]\label{a.dp.def} Let $k$ be a field.
    I follow the terminology of \cite[8.1.12]{MR2964027},  according to which a
     {\it Del~Pezzo surface}
  is a projective $k$-surface $S$ 
  such that $-K_S$ is ample, and  $S_{\bar k}$  is irreducible with   at worst {\it Du~Val   singularites.}

  Du~Val   singularites are also called rational double points or ADE   singularites. 
  In our discussions only  the  $A_n$ singularities $(xy=z^{n+1})$ will appear.

  The {\it degree} of $S$ is the self-intersection number of $-K_S$.

  A  {\it weak  Del~Pezzo surface} is a projective surface $S$  such that  $(-K_S\cdot C)\geq 0$ for every curve,  the self-intersection number of $-K_S$ is positive, and $S_{\bar k}$  has  at worst  Du~Val   singularites.

  For such an $S$  there is a unique
  birational morphism $m_S:S\to S^*$ to a  Del~Pezzo surface $S^*$ such that $K_S=m_S^*K_{S^*}$.
  $S$ and $S^*$ have the same degrees, and
   $m_S$ contracts all the curves $C\subset S$ for which
  $(-K_S\cdot C)= 0$.  If $S$ is smooth, these are the smooth rational curves with self-intersection $-2$.
  $S^*$ is singular at the images of these curves.

  {\it Note on terminology.} Many authors require (sometimes implicitly)  that a Del~Pezzo surface be smooth. For us the more encompassing definition is more convenient. In the construction of the Geiser involutions, the
  most important surface is a singular Del~Pezzo surface; see Example~\ref{a6dp2.exmp}.

  A normal, projective surface  $S$ with $-K_S$ ample
  is called a {\it log Del~Pezzo surface.} Some authors assume that the singularities are quotient singularities, some allow a boundary,  and some drop the modifier `log.'
  So {\it caveat lector.}
 \end{say}

 \begin{defn}[Split nodes] \label{node.defn} Let $C$ be a curve over a field $k$ and
   $c\in C(k)$ a point.  $C$ has a {\it node} at $c$ if,  after completion, it is isomorphic to $\spec  k[[x,y]]/(q)$ where $q$ is a homogeneous degree 2 polynomial with 2 simple roots.

   The node is {\it split} if  $q$ has 2 roots in $k$. Then we can use the standard form $\spec  k[[x,y]]/(xy)$. 
   The node is {\it unsplit} if  the roots of $q$ are conjugate over $k$.
   If $\chr k\neq 2$ then we can use the normal form
   $\spec  k[[x,y]]/(x^2-ay^2)$, where $a$ is not a square in $k$.
    \end{defn}

 \begin{defn}\label{inf.d.defn} Let $D_{\infty}$ denote the {\it infinite dihedral group.}
   It can be given as  $\langle a,b: a^2=1, aba=b^{-1}\rangle$,  as
   $\langle a,c: a^2=c^2=1\rangle$  (setting $c=ab$), as  $\z/2\ast \z/2$, or as
   $\z\rtimes \z/2$ with the $\z/2$-action $r\to -r$  on $\z$. 
 \end{defn}

  \begin{say}[Moduli of log K3 pairs]\label{moduli.lk3.say}
      Let  $S$ be a  Del~Pezzo surface and $N\in |-K_S|$  a geometrically irreducible, nodal curve.
     
      If $S=\p^2$ and $N$ has a split node, then, up to isomorphism $\bigl(\p^2, (xyz=x^3+ay^3)\bigr)$ are the only ones over any field $k$.  If $a$ is a cube in $k$ then we get $\bigl(\p^2, (xyz=x^3+y^3)\bigr)$.

 For $d<9$ the nodal  pairs $(S, N_S)$ with $\deg S=(N_S^2)=d$ form a
$(9-d)$-dimensional family, see \cite{MR3314827} for a detailed treatment of their moduli theory.

A naive argument is the following.
      If $S$ is smooth then we can write $\pi:S\to \p^2$ as the blow-up of $\p^2$ at $c$ points. Each $\pi$-exceptional curve meets $N_S$ with multiplicity 1, so
      $(S, N_S)$ is the blow-up of  $(\p^2, N )$ at  $c$ smooth points on $N $.
      Thus the moduli space is roughly the same as the space of $c$   smooth points in $N$.
      
      Note that if $c=1$ or $c=2$, then $\pi:S\to \p^2$ is unique, but not
for $c\geq 3$.  We have 
finitely many  choices for $\pi$ for $c\leq 8$, and at most countably many for $c\geq 9$.

Note that the Torelli theorem of  \cite[\S 5]{MR3314827}
suggests that the `best' choice may be repeatedly blowing up the
point  $(1{:}{-}1{:} 0)\in (xyz=x^3+y^3)\bigr)$.
This gives the  unique pair where the  mixed Hodge structure on $H_2(S\setminus N_S)$  is split.  However, this is only  a weak Del~Pezzo surface for $8\geq c\geq 2$.
  \end{say}

  We need a few results on discrepancies.
     For general introductions see
     \cite[Sec.2.3]{km-book} or \cite[Sec.2.1]{kk-singbook}.
     However, we  need only 3 simple facts, which are easy to state
     (and prove)  for smooth surfaces.  

  \begin{say}[Discrepancies for surface pairs]\label{discrep.say}
    Let $S$ be a smooth surface and $\pi:T\to S$ a proper, birational morphism.
We are interested in the local behavior of $\pi$ over a point $s\in S$, 
  so   assume that $\pi$ is an isomorphism over $S\setminus\{s\}$, and $T$ is normal.

    We can write $K_T\sim \pi^*K_S+E$ where $E$ is effective and $\pi$-exceptional.

    Next let $C$ be a (possibly reducible) curve on $S$ and $C_T$ its birational transform on $T$. Then  $C_T\sim \pi^*C-D$ where  $D$ is effective and $\pi$-exceptional.  For  any (say rational) $c>0$ we can formally write
    $$
    K_T+cC_T\simq \pi^*(K_S+cC)+E-cD.
    \eqno{(\ref{discrep.say}.1)}
    $$
    The coefficients in $E-cD$ are called the {\it discrepancies} of
    $(S, cC)$.

   The claims below can be established by checking them for 1 blow-up and using induction.
    For precise references  see \cite[2.29]{kk-singbook}
    for (\ref{discrep.say}.2--3) and \cite[2.31]{km-book} or \cite[2.7]{kk-singbook} for
    (\ref{discrep.say}.4).

\medskip
    (\ref{discrep.say}.2)  $E-cD$ is effective for every $\pi$ iff
    $c\cdot \mult_s C\leq 1$. (Such pairs   $(S, cC)$ are called {\it canonical.})

    (\ref{discrep.say}.3)  $E-cD$ is effective  and $\supp(E-cD)=\ex(\pi)$ for every $\pi$ iff
    $c\cdot \mult_s C< 1$. (Such pairs   $(S, cC)$ are called {\it terminal.})

    (\ref{discrep.say}.4)  If $C$ has a node at $s$ and $c\leq 1$ then
all coefficients in $E-cD$ are  $\geq -1$.
(The latter is the {\it log canonical} property.)
    \end{say}

\section{Easy cases of Theorem~\ref{main.thm.1}}

\begin{say}[Preliminary remarks for the proof of Theorem~\ref{main.thm.1}]
  \label{prelim.say}
  Let  $(S, N)$ be a log K3 surface with  irreducible boundary.
  Then
  $\cl (S\setminus N)\cong \cl(S)/\z[-K_S]$.
  In particular, if  $k$ is 
  algebraically closed and $S$ is smooth,
  then the class group of $S\setminus N$ is
  \begin{enumerate}
  \item $\z/3$ if $S\cong \p^2$,
  \item $\z+\z/2$ if $S\cong \p^1\times \p^1$, and
  \item  $\z^{9-d}$ in all other cases, where $d=\deg S$.
  \end{enumerate}
  In particular, if  $\phi:S_1\setminus N_1\cong S_2\setminus N_2$,
  then $\deg S_1=\deg S_2$.

  Let $C\in |-m_1K_{S_1}|$ be an irreducible curve. Then the restriction of
  $C$ to $S_1\setminus N_1$ is linearly equivalent to $0$.
  Thus  the restriction of $\phi_*(C)$
  to $S_2\setminus N_2$ is also linearly equivalent to $0$, hence
  $  \phi_*(C)\in |-m_2K_{S_2}|$ for some $m_2$.
  \smallskip

  \end{say}

\begin{say}[Graphs of birational maps]\label{graphs.say}
 Let  $(S_i, N_i)$ be log K3 surfaces and
 $\phi^\circ: S_1\setminus N_1 \cong S_2\setminus N_2$ an isomorphism that 
 does not extend to an isomorphism
      $\phi:(S_1, N_1)\cong (S_2, N_2)$.
 Let $T$ be the normalization of the closure of the graph of 
 $\phi^\circ$ with projections
 $p_i:T\to S_i$.

 Let  $E_i\subset T$ be the reduced exceptional curve of $p_i$.
 Then $E_i$ is not contracted by $p_{3-i}$, so
 $$
 p_{3-i}: E_i\map N_{3-i} \qtq{is birational.}
 $$
 So  $E_i$ is the birational transform of $N_{r-i}$ on $T$, let us thus write
 $N'_{3-i}:=E_i$.
 
Next we use (\ref{discrep.say}.4), thus
 $$
 \begin{array}{ccl}
   K_T+N'_1+N'_2-D_1& \sim&  p_1^*(K_{S_1}+N_1), \qtq{and}\\ 
   K_T+N'_1+N'_2-D_2& \sim&  p_2^*(K_{S_2}+N_2),
 \end{array}
  \eqno{(\ref{graphs.say}.1)}
$$
 where the $D_i$ are effective and $\supp D_i\subset\supp N'_{3-i}$.
 Since $K_{S_i}+N_i\sim 0$, 
 we conclude that
 $D_1\sim D_2$. Thus in fact  $D_1=D_2=0$, hence
 $$
 p_1^*(K_{S_1}+N_1) \sim K_T+N'_1+N'_2\sim p_2^*(K_{S_2}+N_2).
 \eqno{(\ref{graphs.say}.2)}
 $$

 Note that even if the $S_i$ are smooth, $T$ is usually singular.
 For example, for $S=\p^2$ and $\phi=\sigma_+\circ\sigma_-$, one can compute that  the singularities of $T$ are  $\a^2/\frac17(1,1)$ and  $\a^2/\frac1{48}(1,41)$. 
\end{say}

\begin{lem}\label{must.node.lem}
  Let  $(S_i, N_i)$ be log K3 surfaces with  geometrically irreducible  boundaries, and
  $\phi^\circ: S_1\setminus N_1 \cong S_2\setminus N_2$ an isomorphism
   that 
 does not extend to an isomorphism
      $\phi:(S_1, N_1)\cong (S_2, N_2)$.
  Then the $N_i$ are nodal and 
  $p_i: T\to S_i$ contracts $N'_{3-i}$ to the node $s_i\in N_i$.
  \end{lem}

Proof.  Set $w_i:=p_i(N'_{3-i})$. If $w_i$ is a smooth point of $N_i$
then  (\ref{discrep.say}.2) applies at $w_i$, which contradicts
(\ref{graphs.say}.2). \qed

\begin{say}[Proof of Theorem~\ref{main.thm.1} for $d\leq 3$]
  \label{main.thm.1.pf.1}
  Let $T^r\to T$ be the minimal resolution. By (\ref{must.node.lem})  it dominates the blow-up   $B_{s_i}S_i\to S_i$.
  Let $\bar N_i\subset B_{s_i}S_i$ denote the  birational transform of $N_i$.
  Then the  self intersection of $\bar N_i$ is  $d-4$, and this
  can only decrease with further blow-ups.
  For $T^r\to S_{3-i}$ there must be at least one exceptional $(-1)$ curve,
   and the only candidate is  the birational transform of  $N_{i}$.
  So $d\geq 3$ and
  if $d=3$ then $B_{s_1}S_1=T=B_{s_2}S_2$. \qed
  \end{say}

\begin{say}[Proof of Theorem~\ref{main.thm.1} in the unsplit case]
  \label{main.thm.1.pf.2}
Assume that $(S_1, N_1)$ has an unsplit node.

As we noted in (\ref{main.thm.1.pf.1}), $T^r\to S_1$ dominates the blow-up   $B_{s_1}S_1\to S_1$. The boundary of $B_{s_1}S_1$ is $F_1+\bar N_1$, where $F_1$ is the exceptional curve.  $F_1+\bar N_1$ has 2 singular points, and these are conjugate if the node is unsplit. Thus the image of $N'_{3-i}$ must be $F_1$.
Thus again we conclude that $B_{s_1}S_1=T$ and $d=3$. \qed
 \end{say}

\section{Geiser involutions}
  
Geiser involutions are constructed using a natural automorphism of
degree 2 Del~Pezzo surfaces. We start by recalling their properties; see 
\cite[Chap.8]{MR2964027} for details.
We assume for simplicity that the characteristic is $\neq 2$.

  \begin{say}[Degree 2 Del~Pezzo surfaces]\label{a.d2dp}
    Let $T$ be a  degree 2 Del~Pezzo surface.
The linear system $|-K_T|$ gives a morphism  $p:T\to \p^2$, which 
    is a double cover, branched along a
  quartic curve $R\subset \p^2$.
  We let $\tau$ or $\tau_T$ denote the covering involution.

  The anticanonical class is the pull-back of the line class $L$ on $\p^2$.

  $T$ is smooth iff $R$ is smooth.
  The singularities of $T$ correspond to the singular points of $R$, and they have the same name in the ADE classification.
  
  A curve $E\subset T$ is called a {\it line} if $(-K_T\cdot E)=1$.
  Since $-K_T=p^*L$, we see that $E$ is a line iff
  $p(E)\subset \p^2$ is a line and $p^{-1}(p(E))=E+E'$ is a sum of 2 lines.
    This happens over the bitangents of $R$. If $R$ is smooth, there are 28 bitangents, giving 56 lines on $T$.
  (If $p(E)\subset R$ then $p^{-1}(p(E))=2E$.)

  If $E,E'$ are 2 lines and $p(E)\neq p(E')$, then $E\cap E'$ is at most 1 point.
  Thus, if $E\cap E'$ consists of 2 points, then $p(E)= p(E')$, hence
  $\tau(E)=E'$.

  Let $D\subset T$ be a curve. Then $D+\tau_T(D)=p^{-1}(p(D))$, hence
  $$
  D+\tau_T(D)\in |-r K_T| \qtq{where}  r=(-K_T\cdot D).
  \eqno{(\ref{a.d2dp}.1)}
  $$
  Let  $T$ be a smooth, weak  Del~Pezzo surface of  degree 2
and  $m_T:T\to T^*$ the corresponding Del~Pezzo surface of  degree 2.
  Then the
  covering involution $\tau^*$  on $T^*$ lifts to an involution $\tau$ on $T$.

  {\it \ref{a.d2dp}.2 Note on  characteristic 2.}
  In the cases that we consider later, $T\to \p^2$ is always separable, so
  the above descriptions work.

   \end{say}

\begin{say}[Geiser involutions]  These were discovered in \cite{geiser}; see \cite[Sec.8.7]{MR2964027} for a modern treatment.

  Start with $S=\p^2$ and blow up $7$ points $P=\{p_1, \dots, p_7\}$ to get  $\pi_P:T_P\to S$.
If no 3 points are on a line and  no 6 on a  conic, then  $T_P$ is a smooth  Del~Pezzo surface of degree 2. 
Thus $T_P$ is a double cover of (another copy of)  $\p^2$, branched along a smooth quartic curve. Let $\tau_P:T_P\to T_P$ be the covering involution of $T_P\to \p^2$.
 The classical
 {\it Geiser involutions} are given as
 $$
 \sigma_P:=\pi_P\circ\tau_P\circ \pi_P^{-1}: \p^2\map \p^2.
 $$
 More generally, let $S$ be a Del Pezzo surface
 and $\pi_T:T\to S$  a birational morphism from a smooth,  degree 2 Del~Pezzo surface to $S$. Let  $\tau_T:T\to T$ be the  covering involution. 
 As above, we get
$$
 \sigma_T:=\pi_T\circ\tau_T\circ \pi_T^{-1}: S\map S.
 $$
  I also call these 
{\it Geiser involutions.}  
(Again, the terminology is  not uniform.
If $S=\p^2$ but $T\neq T^*$, then  \cite{Bayle-Beauville00} considers $\sigma_T$ to  be a
{\it de~Jonqui\`eres involution.}
For the current purposes, the representation $\sigma_T=\pi_T\circ\tau_T\circ \pi_T^{-1}$  shows the underlying geometry much more clearly.)
\end{say}

Next we look at a special case where the blow-up centers of $\pi_T$ are in a rather degenerate position.

\begin{say}[Geiser involutions from split  nodes] \label{a.obu}
  Assume now that  $S$ is a Del~Pezzo surface of degree $r+2\geq 3$, and $N\in |-K_S|$ a geometrically irreducible curve with a split node at a smooth point $s\in S$.

  Pick a local branch  $s\in B_+\subset N$ and
  let $\pi_+:T_+\to S$ denote the blow-up of $s\in B_+$ successively
  $r$ times. We check below that $T_+$ is a weak Del~Pezzo surface of degree 2.
  Thus we get a Geiser involution   $\sigma_+ :S\map S$.
  
  Starting with the other branch $B_-\subset N$, the first blow-up is the same, but the next $r-1$ blow-ups are different. For $r\geq 2$ this gives another  Geiser involution
  $\sigma_- :S\map S$. (If $r=1$ then we blow up $s$ only once, hence the
 2 involutions coincide.)

  As we noted after  (\ref{gei.i.i}.1), the branches $B_{\pm}$ are
  indistinguishable from each other,  so
  one should think of $\{\sigma_+, \sigma_-\}$ as a pair of involutions.

  For future computations, let us give a detailed description of
   $\pi_{\pm}:T_{\pm}\to S$ and the action of the covering involution $\tau$.

  More generally, let $S$ be a  surface and  $N\subset S$ a  curve with a split node $s\in N$.
Assume that $S$ is smooth at $s$.
  We represent this as
  $B {\mbox{ --- }} B$ where the $B$ are the local branches at the node. In suitable local analytic coordinates we can identify the branches with the coordinate axes.
  This is very covenient for explicit computations.

 All diagrams below are arranged so that we blow-up points on the  branch $B$ that is on the right hand side.

For example,  3 blow-ups give the following sequence of dual graphs, where we denote an exceptional curve by $E$,  and the subscript  tells the number of blow-ups needed for that curve to appear. 
$$
B {\mbox{ --- }} B
\qquad
B {\mbox{ --- }} E_1  {\mbox{ --- }}B
\qquad
B {\mbox{ --- }} E_1  {\mbox{ --- }} E_2  {\mbox{ --- }}B
\qquad
B {\mbox{ --- }} E_1  {\mbox{ --- }}E_2  {\mbox{ --- }} E_3  {\mbox{ --- }}B
$$
After $r$ blow-ups a more symmetric version is 
$$
\begin{array}{ccccc}
  E_1 & - &E_2\ -\cdots -\ E_{r-2}& - & E_{r-1} \\
  \vert &&&& \vert\\
  N_T && \rule[1pt]{40pt}{0.4pt} && {\ }E_r
\end{array}
\eqno{(\ref{a.obu}.1)}
$$
where $N_T$ denotes the birational transform of $N$.

Since we repeatedly blow up nodes,
  the sum of all curves in (\ref{a.obu}.1) is a member of $|-K_T|$.
  Thus $T$ is a weak Del~Pezzo surface of degree 2.
  Note that $(-K_T\cdot E_i)=0$ for $1\leq i<r$,  $(-K_T\cdot E_r)=1$
  and $(-K_T\cdot N_T)=1$.

  Let $m_T:T\to T^*$ be the
  birational morphism  to a degree 2 Del~Pezzo surface $T^*$ as in
  (\ref{a.d2dp}). Then
  \begin{enumerate}\setcounter{enumi}{1}
    \item $m_T$ contracts the curves $E_1, \dots E_{r-1}$,
\item   $m_T(N_T)$ and $m_T(E_r)$ are lines that meet at 2 points, and so
  \item  $\tau^*$ interchanges $m_T(N_T)$ and $m_T(E_r)$.
\item   Therefore $\tau(N_T)=E_r, \tau(E_r)=N_T$, and $\tau(E_i)=E_{r-i}$ for $i=1,\dots, r-1$.
  \end{enumerate}

  {\it  \ref{a.obu}.6 Note on  characteristic 2.}
  Since $\tau^*$ interchanges $m_T(N_T)$ and $m_T(E_r)$, the anticanonical map
  $T\to \p^2$ is separable.
  
\end{say}

\begin{exmp}\label{a6dp2.exmp} If $S=\p^2$ then $r=7$, hence 
  $m_T:T\to T^*$ contracts the curves $E_1 - \cdots - E_6$.
  Thus $T^*$ has an $A_6$ singularity.
  The branch curve of $T^*\to \p^2$ is then a 
  plane quartic with an $A_6$ singularity.
There is a unique such curve up to isomorphism, which
can be given by an affine equation  $(y-x^2)^2-xy^3=0$.
The  corresponding affine equation of $T^*$ is $u^2=(y-x^2)^2-xy^3$. 

There are 4 lines on $T^*$ that pass through the singular point.
$L_{\pm}$ lie over $(x=0)$ and $M_{\pm}$ lie over $(y=0)$.
Note that $M_+$ and $M_-$ intersect only over the origin $(x=y=0)$, while
$L_+$ and $L_-$ intersect  over  $(x=y=0)$ and  at infinity.
Thus the $L_{\pm}$ are the images of $N_T+E_7$.

Note also that, if $k$ contains a 3rd root of unity  $\epsilon\neq 1$, then $T$ and $T^*$ have an order 3 automorphism
$
(x,y, u)\mapsto (\epsilon x, \epsilon^2 y,  \epsilon^2 u).
       $
\end{exmp}

\begin{say}[Geiser involutions and Yoshihara twists] \label{a.giot}
 A version of the Geiser involution (\ref{a.obu}) for $\p^2$ was considered in
 \cite{MR0769779}, rediscovered  in \cite{MR1942244} and then
 generalized to other Del~Pezzo surfaces in \cite{mcduff2024singular}.

 The constructions in these papers 
  start with the same $r$ blow-ups  as in (\ref{a.obu}), and
 then contract the curves
 $N_T, E_1, \dots, E_{r-1}$, starting with $N_T$.

 The end result is then a birational map  $\phi:(S, N)\map (S', N')$,
 where  $S'$ is a Del~Pezzo surface and $N'$ is the image of $E_r$.
 It is clear that $\deg S=\deg S'$ and that 
 $N'$ is also a  nodal cubic.

 For $S=\p^2$ there is a  unique  pair  $(\p^2, N)$ up to isomorphism over $\c$, hence then $(\p^2, N')\cong (\p^2, N)$.
There are 6 distinct 
isomorphisms  $\rho:(\p^2, N')\cong (\p^2, N)$.
The maps in \cite{MR0769779, MR1942244} are  the Geiser involutions composed with an
automorphism of  $(\p^2, N)$ that switches the 2 branches at the node,  giving  an infinite order automorphism  of $\p^2\setminus N$.

For $d\in\{5,6,7\}$ there is a unique
Del~Pezzo surface of degree $d$ over $\c$, up to isomorphism. Thus, in these cases, 
  $S'\cong S$, and it is easy to see that the same holds for  $d=8$.
However, for $d\in\{5,6,7,8\}$
 the pairs  $(S, N)$ have moduli  (\ref{moduli.lk3.say}), and it is not clear from the descriptions in \cite{MR0769779, MR1942244, mcduff2024singular}  that
 $(S, N)\cong (S', N')$.
 
 The main advantage of working with  Geiser involutions  is  that we get  automorphisms of $S\setminus N$, and the construction works over arbitrary fields.
 
 Note that \cite{mcduff2024singular} studies  deformation invariant properties of the pairs  $(S, N)$, so knowing that  $(S, N)\cong (S', N')$  is not important for their purposes.
\end{say}

\begin{say}[Automorphisms of log K3 surfaces]\label{nodal.p.say}
  The automorphism group of $(S, N)$  acts on $N$, and $\aut(N)\cong \gm\rtimes \z/2$.  The action also fixes $\o_S(K_S)|_N\in \pic^d(N)$ where  $d=\deg S$.
  Note that $\lambda\in \gm=\aut^\circ(N)$ acts on $ \pic^d(N)$ by multiplication by $\lambda^d$.  Thus, for $d\neq 0$,  the subgroup of $\aut(N)$ that leaves 
  $\o_S(K_S)|_N\in \pic^d(N)$  invariant is the dihedral group
   $D_{2d}\cong \z/d\rtimes \z/2$.
Thus we get a homomorphism
$\aut(S, N)\to D_{2d}$.
If $d\geq 1$ then the kernel is cyclic.

If $k$ contains a 3rd root of unity  $\epsilon\neq 1$, then  $\aut\bigl(\p^2, (xyz=x^3+y^3)\bigr)\cong S_3$, generated by the involution $\rho_2:(x:y:z)\mapsto (y:x:z)$ and
 the order 3 element  $\rho_3:(x:y:z)\mapsto (\epsilon x: \epsilon^2y:z)$.  
 As another example, if $\chr k\neq 2$ then $\aut\bigl(\p^1\times \p^1, (x_0^2y_1^2+x_1^2y_0^2=x_0^2y_0^2)\bigr)$ is the  dihedral group $D_8$.  
 For $d\leq 8$ the general pair  $(S, N)$ has no automorphisms.

 Let $T_{\pm}\to S$ be the blow-ups involved in the construction of the Geiser involutions. If   $\rho\in  \aut(S,N)$
 interchanges the branches $B_{\pm}$, then  $\rho$ lifts to an isomorphism
 $\bar\rho: T_{+}\cong T_{-}$. Otherwise
 $\rho$ lifts to an isomorphism $\bar\rho: T_{+}\cong T_{+}$.

Therefore, if $\rho\in  \aut(S,N)$ interchanges the branches $B_{\pm}$, then
$\rho^{-1}\sigma_{\pm}\rho=\sigma_{\mp}$,
and $\rho^{-1}\sigma_{\pm}\rho=\sigma_{\pm}$ otherwise.
This gives the following.
\end{say}

\begin{cor}\label{sd.prod.cor}
Let  $S$ be a Del~Pezzo surface of degree $\geq 4$, and $N\in |-K_S|$ with a split node at a smooth point $s\in S$.
  Then the  subgroup of 
  $\aut(S\setminus N)$ generated by the $\sigma_{\pm}$ and $\aut(S,N)$ is a
  semidirect product
  $\langle \sigma_+, \sigma_-\rangle\rtimes \aut(S,N)$. \qed
\end{cor}

       Note that $\aut(S, N)$ can be  infinite if $(N^2)\leq 0$; see   \cite{li2022cone} for a description of these groups and applications to the Morrison cone conjecture
       for log K3 surfaces.

\section{Cusps}

      \begin{defn}[Cusp]\label{cusp.defn}  Let $C$ be  
        a geometrically reduced curve  over a field $k$ with normalization  $\pi:\bar C\to C$. For a point $c\in C(k)$ set  $\bar c:=\red \pi^{-1}(c)$.
        We say that $C$ has 
        a {\it cusp} or is {\it unibranch}  at $c\in C$ if 
        $k(\bar c)=k(c)=k$.

       Let $S$ be a smooth  surface germ and fix local  (analytic or formal)
           coordinates $(x,y)$ at $s\in S$. Then   a cusp   $s\in C\subset S$
        can be parametrized as
        $$
        t\mapsto \bigl(x(t):=t^p+(\mbox{higher terms}), y(t):=t^q+(\mbox{higher terms})\bigr),
        \eqno{(\ref{cusp.defn}.1)}
        $$
        (after scaling $x, y$ if necessary).
        The values $p, q$ are  the {\it contact orders} of $C$ with the axes:
       $(C\cdot (\mbox{$x$-axis}))_s=q$ and
       $(C\cdot (\mbox{$y$-axis}))_s=p$.
        
        The multiplicity of $C$ is $\min\{p,q\}$, independent of the choice of the coordinates.  If $q\nmid p$  and $p\nmid q$ then both $p, q$ are invariants of $C$;  we then call $C$ a   {\it $(p,q)$-cusp.}
        
        Note that if $p=rq$ then the coordinate change $x':=x- y^r$
        gives $\mult x'(t)>\mult x(t)$.
        Newton proved that if $C$ is singular, then after finitely many steps we reach a
        $$
    t\mapsto \bigl(x(t):=t^{p'}+(\mbox{higher terms}), y(t):=t^{q}+(\mbox{higher terms})\bigr),
    $$
    such that  $q\nmid p'$;
    see \cite[Sec.8.3]{br-kn} or  \cite[Sec.1.1]{k-res}.
    Then  $C$ is a $(p',q)$-cusp.
    
   After blowing up the origin,   the parametrization (\ref{cusp.defn}.1) changes to 
   $$
   \begin{array}{lcll}
     (x/y,y)&=&\bigl(t^{p-q}+(\mbox{higher}), t^q+(\mbox{higher})\bigr) &\mbox{if $p>q$, and}\\[1ex]
   (x,y/x)&=&\bigl(t^p+(\mbox{higher}), t^{q-p}+(\mbox{higher})\bigr)
     &\mbox{if $p<q$.}\end{array}
   \eqno{(\ref{cusp.defn}.2)}
   $$
  \end{defn}

      \begin{say}[Blowing up  cusps] \label{a.catb} Let $s\in N\subset S$ be a nodal curve; call one branch $B_+$ and the other  $B_-$.
        Using these branches as local coordinate axes, 
   for a cusp  $s\in C\subset S$ the key invariants are the
    {\it contact orders}  with $N$. We write these as 
    $p:=(B_-\cdot C)_s$ and $q:=(B_+\cdot C)_s$, and 
     represent this with the  dual graph  $B_- \stackrel{p,q}{{\mbox{ --- }}} B_+ $.

 After blowing up $s$, the  dual graphs change as 
$$
\circ \stackrel{p,q}{{\mbox{ --- }}} \circ
\qquad\rightsquigarrow\qquad
\begin{array}{ll}
  \circ \stackrel{p-q,q}{{\mbox{ --- }}} E_1 {\mbox{ --- }}\circ &\mbox{if $p\geq q$, and}\\[1ex]
 \circ {\mbox{ --- }} E_1 \stackrel{p,q-p}{{\mbox{ --- }}}\circ& \mbox{if $p\leq q$.}
\end{array}
\eqno{(\ref{a.catb}.1)}
$$
Note that if $p=q$ then the 2 versions mean the same thing: the transform of $C$ meets $E_1$  away  from the transform of $N$.

By induction, after $r$ blow-ups of $s\in B_+$ we get the following dual graphs.
  $$
  \begin{array}{ll}
 (\ref{a.catb}.2)\quad    B \stackrel{p-q,q}{{\mbox{ --- }}} E_1 {-} \cdots {-} E_r {\mbox{ --- }} B
    &
  \mbox{if $q\leq p$,}\\[1ex]
  (\ref{a.catb}.3)\quad   B  {-} \cdots {-} E_i \stackrel{(i+1)p-q, q-ip}{{\mbox{ --- }}} E_{i+1} {-} \cdots {-}  B
  &
  \mbox{if $p\leq ip\leq q< (i+1)p\leq rp$,}\\[1ex]
   (\ref{a.catb}.4)\quad  B {\mbox{ --- }} E_1 {-} \cdots {-} E_r \stackrel{p,q-rp}{{\mbox{ --- }}} B
    &
   \mbox{if $rp\leq q$.}
   \end{array}
   $$

\end{say}

\section{Geiser involutions and  cusps}

\begin{prop} \label{a.catbt}
  Let $S$ be a   Del~Pezzo surface of degree $r+2\geq 4$. Let $N\in |-K_S|$ be a geometrically irreducible  curve with a split node $s\in N$ such that  $S$ is  smooth at  $s$,
  and  $\sigma_{\pm}$ the corresponding Geiser involutions  (\ref{a.obu}).
  Let $s\in C^\circ\subset S$ be a (germ of a) cusp represented with the  dual graph
   $B \stackrel{p,q}{{\mbox{ --- }}} B$. 
   Then  $\sigma_+(C^\circ)$ is  represented with the  dual graph
  $$
    \begin{array}{ll}
      B \stackrel{p, rp-q}{{\mbox{ --- }}} B
       &
      \mbox{if $q\leq rp$, and}\\[1ex]
       B \stackrel{q-rp,p+r(q-rp)}{{\mbox{ --- }}} B
       &
      \mbox{if $q\geq rp$.}
    \end{array}
    \eqno{(\ref{a.catbt}.1)}
    $$
\end{prop}
Note that if $q=rp$ then both forms mean that intersection point is different from $s$.  Since $N\setminus\{s\}$ is connected, it does not make sense to ask
on which branch of $N$  the point is.

\medskip

Proof. After $r$ blow-ups we have the dual graphs (\ref{a.catb}.2--4).
The action of  $\tau_B$ is described in   (\ref{a.obu}.2--5). 
Thus  the cases  (\ref{a.catb}.2--4) are transformed into the following.
$$
  \begin{array}{ll}
 (\ref{a.catbt}.2)\quad    B  {-} \cdots {-} E_{r-1} \stackrel{q,p-q}{\mbox{ --- }} E_r {\mbox{ --- }} B
    &
  \mbox{if $q<p$,}\\[1ex]
(\ref{a.catbt}.3)\quad   B  {-} \cdots {-} E_{r-i-1} \stackrel{q-ip, (i+1)p-q}{{\mbox{ --- }}} E_{r-i} {-} \cdots {-}  B
  &
  \mbox{if $p\leq ip<q< (i+1)p\leq rp$,}\\[1ex]
  (\ref{a.catbt}.4)\quad    B {\mbox{ --- }} E_1 {-} \cdots {-} E_r \stackrel{q-rp,p}{{\mbox{ --- }}} B
    &
   \mbox{if $rp<q$.}
   \end{array}
  $$
  Next we contract the curves $E_r, \dots, E_1$.
In cases (\ref{a.catbt}.2--3)
  the first contractions that involve the curve $C$ are 
  $$
   \begin{array}{lcll}
  E_{r-1} \stackrel{q,p-q}{\mbox{ --- }} E_r 
&\rightsquigarrow&
  E_{r-1} \stackrel{p,p-q}{\mbox{ --- }} B
  &\mbox{and}\\[1ex]
  E_{r-i-1} \stackrel{q-ip, (i+1)p-q}{{\mbox{ --- }}} E_{r-i}
  &\rightsquigarrow&
  E_{r-i-1} \stackrel{p,(i+1)p-q}{{\mbox{ --- }}} B. & 
  \end{array}
  $$
   These correspond to the first case in (\ref{a.catb}.1).
   After that all contractions are as in the second case in (\ref{a.catb}.1).
   Thus $p$ stays fixed and we add $(r-i-1)p$ to the second  component, to end with
   $(p, rp-q)$.
   
   In case (\ref{a.catb}.1) we add $r(q-rp)$ to $p$. 
   \qed

\begin{cor}\label{a.thm}
  Let $S$ be a   Del~Pezzo surface of degree $r+2\geq 4$. Let $N\in |-K_S|$ be a geometrically irreducible  curve with a split node $s\in N$ such that  $S$ is  smooth at  $s$.
  Let $\sigma_+, \sigma_-$ the corresponding Geiser involutions.
 
    Let $C^\circ\subset S$ be a germ of a  cusp  that has  contact orders $(p, q)$ with the branches of  $N$. Then $\sigma_+(C^\circ)$ and $\sigma_-(C^\circ)$ have  contact orders 
\begin{enumerate}
\item   $(p, rp-q)$  and $ (rq-p,q)$  if $\tfrac1{r}<\tfrac{p}{q}<r$, 
\item    $\bigl(q-rp, p+r(q-rp)\bigr)$ and $ (rq-p,q)$ if
  $\tfrac{p}{q}\leq \tfrac1{r}$, 
 \item    $(p, rp-q)$ and $\bigl(q+r(p-rq), p-rq\bigr)$
  if $r\leq \tfrac{p}{q}$. \hfill \qed
\end{enumerate}
  \end{cor}

\begin{say}[Proof of Lemma~\ref{cusp.transf.cor}]\label{cusp.transf.cor.pf}
  Using the notation of Lemma~\ref{cusp.transf.cor}, let $C\subset S$ be an irreducible curve with only 1 point on $N$ having contact orders  $(p,q)$ at $s$.
  Then  $\mult(C)=\min\{p,q\}$, and
  $\ndeg (C)=\tfrac1{d}(-K_S\cdot C)$.
  The rest follows by rewriting the formulas (\ref{a.thm}.1--3) in terms of
  $\mult(C)$ and $\ndeg (C)$. \qed
\end{say}

  \section{Upper bound for the multiplicity}\label{up.bd.sec}

  We check that the upper bound $\mult(C)<2 \cdot\ndeg (C)$
  needed in 
(\ref{bounds.i.say}.3) holds for almost all curves.

    \begin{say}[Max Noether's genus formula]\label{genus.say}
(See \cite[Sec.9.2]{br-kn} or \cite[p.38]{k-res}.)
      The {\it genus} of a  planar curve singularity is  $\delta(0\in C):=\sum \frac{m_i(m_i-1)}{2}$,  where $m_i$ are the multiplicities of all points that appear after some blow-ups (called infinitely near points).
    If $C$  is an irreducible, projective  curve on a smooth surface $S$  with normalization $\bar C$, then Noether's formula is
    $$
    g(\bar C)=\tfrac12\bigl(C\cdot (C+K_S)\bigr)+1-\tsum_{c\in C} \delta(c\in C).
    \eqno{(\ref{genus.say}.1)}
    $$
    Let $S$ be a smooth, weak Del~Pezzo surface of degree $d$ and
    $C\subset S$ an irreducible curve.  Write
    $C\simq -mK_S+F$ where $(K_S\cdot F)=0$.
    Then
    $$
    \bigl(C\cdot (C+K_S)\bigr)= dm(m-1)+(F^2). 
    $$
Thus we get that 
$$
2\cdot \tsum_{c\in C} \delta(c\in C)=dm(m-1)+2-2g(\bar C)+(F^2).
     \eqno{(\ref{genus.say}.2)}
    $$      
Using (\ref{cusp.defn}.2) and  induction we see  that
    $$
    2\cdot \delta\bigl(\mbox{$(p, q)$-cusp}\bigr)\geq (p-1)(q-1),
    \eqno{(\ref{genus.say}.3)}
    $$
    and equality holds iff $(p,q)=1$. Putting these together gives the following.

    \medskip

    {\it Claim \ref{genus.say}.4.} Let $S$ be a smooth, weak Del~Pezzo surface of degree $d$ and
    $C\subset S$ an irreducible curve such that
    $C\simq -mK_S+F$ where $(K_S\cdot F)=0$. If $C$ has a $(p, q)$-cusp then
    $$
    (p-1)(q-1)\leq dm(m-1)+2,
    $$
    with equality holding iff $(p,q)=1$, $C$ is rational with  at most 1 singularity, and   $C\simq -mK_S$ for some $m$. \qed
        \end{say}

     \begin{cor}  \label{q.2m.cor} Let $S$ be a smooth, weak Del~Pezzo surface of degree $d\geq 4$, and
       $C\subset S$ an irreducible curve with a cusp at $\{s\}=C\cap N$.  Then one of the following holds.
       \begin{enumerate}
       \item $\mult_s C< \frac2{d}(-K_S\cdot C)$, or 
       \item  $d\in\{4,5\}$,  and $C\in |-K_S|$ has a $(2,3)$ cusp at $s$.
               \end{enumerate}
  \end{cor}

     Proof.  If $C$ and $N$ have contact orders  $(p,q)$ at $s$, then 
      $p+q=dm$  where $m=\frac1{d}(-K_S\cdot C)$, and  $(p-1)(q-1)\leq dm(m-1)+2$ by (\ref{genus.say}.4). Thus
  $$
  (dm-q-1)(q-1)\leq dm(m-1)+2.
  \eqno{(\ref{q.2m.cor}.3)}
  $$
  The left hand side reaches its maximum at $  q_0=\frac12(dm+2)$, and is monotone decreasing as a function of $|q-q_0|$. So it is enough to check what happens at $q=2m$. Then
  $$
  (dm-2m-1)(2m-1)=2dm^2-4m^2-dm+1=dm(m-1)+1+(d-4)m^2.
  $$
  This is bigger than $dm(m-1)+2$, except when  $d=5$ and $m=1$, or when  $d=4$.

  If $d=4$ then $p+q=4m$ gives (\ref{q.2m.cor}.1) except when $p=q=2m$.
  Then we have a $(2m, 2m+1)$ cusp or worse. (\ref{genus.say}.4) then becomes
  $(2m-1)2m\leq 4m(m-1)+2$ which implies that $m=1$.
  \qed

  \section{Noether-Fano method}\label{n-f.meth.sec}

Here we study curves with a low multiplicity cusp at the node of $N$, closely following \cite[Sec.5.1]{ksc}.

\begin{thm}\label{n-f.meth.thm}
  Let $S_i$ be Del~Pezzo surfaces, $s_i\in N_i\in  |-K_{S_i}|$ geometrically irreducible,  nodal curves and  $\phi^\circ:(S_1\setminus N_1)\cong (S_2\setminus N_2)$  an isomorphism.
  Let $N_1\neq C_1\subset S_1$ be an irreducible curve
  and  set $C_2:=\phi^\circ_*(C_1)$. Assume that
  $$
  \mult_{s_i}C_i\leq \tfrac1{\deg S_i}(-K_{S_i}\cdot C_i)\qtq{for $i=1,2$.}
  \eqno{(\ref{n-f.meth.thm}.1)}
  $$
  Then
  \begin{enumerate}\setcounter{enumi}{1}
  \item  $( K_{S_1}\cdot C_1)=(K_{S_2}\cdot C_2)$, and
    \item if one of the inequalities in (\ref{n-f.meth.thm}.1) is strict, then
      $\phi^\circ$ extends to an isomorphism  $\phi:S_1\cong S_2$.
  \end{enumerate}
\end{thm}

Proof.   Set $m_i:=\frac1{\deg S_i} (-K_{S_i}\cdot C_i)$. We may assume that $m_1\leq m_2$.

Using the notation in (\ref{graphs.say}), 
let $p_i:T\to S_i$ be the projections and $C_T$ the birational transform of the $C_i$.  Write
$C_i\simq  -m_iK_{S_i}+F_i$.

$K_T\simq p_i^*K_{S_i}+E_i$ and $C_T=p_i^*C_i-D_i$.
As in (\ref{discrep.say}.1) for any $c$  we have
 $$
 K_T+cC_T\simq    p_i^*\bigl(K_{S_i}+c(C_i-F_i)\bigr)+E_i-cD_i +cp_i^*F_i.
 $$
  Now choose  $c=\frac1{m_1}$. Then
 $$
 K_{S_1}+c(C_1-F_1) \simq 0
 \qtq{and}
K_{S_2}+c(C_2-F_2)\simq \tfrac{m_1-m_2}{m_1}K_{S_2}.
$$
Hence
 $$
 p_2^*\bigl(\tfrac{m_1-m_2}{m_1}K_{S_2}\bigr)+E_2-cD_2 +cp_2^*F_2
 \simq E_1-cD_1 +cp_1^*F_1.
 $$
 Take intersection number with  $-p_2^*K_{S_2}$. We get that
 $$
 \tfrac{m_1-m_2}{m_1}\deg S_2=\bigl(-K_{S_2}\cdot (p_2)_*(E_1-cD_1)\bigr).
 $$
 Assumption (\ref{n-f.meth.thm}.1) and (\ref{discrep.say}.2) imply that  $E_1-cD_1$ is effective.
 Therefore 
 $m_1\geq m_2$, hence $m_2=m_1$.

 If strict inequality holds in (\ref{n-f.meth.thm}.1) then
 $(p_2)_*(E_1-cD_1)$ is a positive multiple of $N_2$ by (\ref{discrep.say}.3), 
proving (\ref{n-f.meth.thm}.3).
\qed

\begin{thm}\label{n-f.meth.thm.cor}
  Let $S_i$ be Del~Pezzo surfaces, $s_i\in N_i\in  |-K_{S_i}|$ geometrically irreducible, nodal curves and  $\phi^\circ:(S_1\setminus N_1)\cong (S_2\setminus N_2)$  an isomorphism.
  Assume that  $N_2$ has a split node.  Then there is a $\psi\in \langle \sigma_+, \sigma_-\rangle\subset \aut(S_2\setminus N_2)$  such that  $\psi\circ \phi^\circ$ extends to an isomorphism  $S_1\cong S_2$.
\end{thm}

Proof. Let $C_1\in |-K_{S_1}|$  be the curve $E$ in (\ref{E.1pt.lem}),
which is defined over $\bar k$.
Set  $C_2:=\phi^\circ_*(C_1)$.

Choose  $\psi\in\langle \sigma_+, \sigma_-\rangle$ such that
$\ndeg\bigl(\psi_*(C_2)\bigr)$ is the smallest possible.
 Since $E$ is not rational, we are in the first case of Corollary~\ref{cusp.mult.thm}, hence 
$$
\mult\bigl(\psi_*(C_2)\bigr)\leq \ndeg\bigl(\psi_*(C_2)\bigr).
\eqno{(\ref{n-f.meth.thm.cor}.1)}
$$
Since $\sigma_+, \sigma_-$ are defined over $k$,
$\psi\circ \phi^\circ$  is defined over $k$.

We claim that $\psi\circ \phi^\circ$ extends to an isomorphism  $S_1\cong S_2$.
A map is an  isomorphism  over $k$ iff it is an  isomorphism  over $\bar k$.
Thus we can now work over $\bar k$ and 
 apply (\ref{n-f.meth.thm}) to  $\psi\circ \phi^\circ$.
Note that $s_1\notin C_1$, so
(\ref{n-f.meth.thm}.1) is a strict inequality for $i=1$.
Next (\ref{n-f.meth.thm}.1)   for $i=2$ is equivalent to
(\ref{n-f.meth.thm.cor}.1).

Thus  $\psi\circ \phi^\circ$ extends to an isomorphism  $S_1\cong S_2$.
\qed

\begin{cor}\label{n-f.meth.thm.cor.2}
  Let $S$ be a Del~Pezzo surface of degre $\geq 4$,  and  $N\in  |-K_{S}|$ a geometrically irreducible, nodal curve with  a split node at a smooth point $s\in S$.
  Then the $\sigma_{\pm}$ and  $\aut(S, N)$ generate
  $\aut(S\setminus  N)$. \qed
\end{cor}

 \begin{lem} \label{E.1pt.lem}
   Let  $(S, N)$ be a log K3 surface with  geometrically irreducible boundary.
   Then, over $\bar k$,   there is an 
   irreducible curve $E\in  |-K_S|$ such that 
      $N\cap E=\{b\}$ is a single  point and $E$ is smooth at $b$.
 \end{lem}

 Proof.  There is a point $b\in N(\bar k)$ such that
 $\o_N(d[b])\cong \o_S(-K_S)|_N$. (If $\chr k\nmid d$, then there are $d^2$ such points if $N$ is smooth and $d$ such points if $N$ is nodal.)  From the sequence
 $$
 0\to \o_S\to  \o_S(-K_S)\to  \o_S(-K_S)|_N\to 0
 $$
 and $H^1(S, \o_S)=0$ we get   a pencil of curves $|E|\subset |-K_S|$ meeting $N$ only at  $b$.
Since $N\in  |E|$ and $N$ is smooth at $b$, general members of $|E|$ are smooth at $b$.
 \qed

 {\it Remark \ref{E.1pt.lem}.1.} Note that $E$ is automatically smooth in characteristics $\neq 2,3$, but a  more careful inspection shows that 
$E$ is actually smooth even in characteristics $ 2,3$.

\section{Affine lines in $(\p^2, N )$}  \label{al.p2.sec}

We aim to describe all embeddings  $\a^1\into (\p^2\setminus N )$.
The closure of the image is then a cuspidal curve $C\subset \p^2$ that
meets $N $ at a single point.

In this section we work over a field of characteristic 0. The main reason for this is the use of
\cite{MR1016092} in (\ref{1/3.1/2.lem}), but there are also other problems in characteristics 2 and 3.

\begin{thm}\label{autom.thm.reduce}
  Let  $N\subset \p^2$ be a geometrically irreducible cubic with a split node, and
  $C\subset \p^2$  a cuspidal, rational curve that meets $N $ at a single point $c$. Assume that all cusps of $C\setminus\{c\}$ have multiplicity $\leq 2$.

  Then there is a $\phi\in \langle \sigma_+, \sigma_-\rangle$ such that
  $\phi_*(C)\subset \p^2$ is  a line or a conic.
\end{thm}

Lines and conics meeting $N $ at a single point are listed in (\ref{3.smooth.lines}).

  The theorem implies that  $C\setminus\{c\}$ is actually smooth.
  It is quite possible that 
  the extra assumption (multiplicity $\leq 2$) is not needed.

  \medskip
  Proof. By Corollary~\ref{cusp.mult.thm} there is a   $\phi\in \langle \sigma_+, \sigma_-\rangle$    such that
  $\mult \bigl(\phi(C)\bigr)\leq  \ndeg\bigl(\phi(C)\bigr)$.
  Note that $\ndeg $ is the normalized degree,
    so for the usual degree  this says that
  $$
  \mult_s \bigl(\phi(C)\bigr)\leq  \tfrac13\deg\bigl(\phi(C)\bigr).
  \eqno{(\ref{autom.thm.reduce}.1)}
  $$
  If $\mult_s \bigl(\phi(C)\bigr)\geq 2$ then
  $s$ is a maximal multiplicity point of $\phi(C)$, and then
  (\ref{autom.thm.reduce}.1) contradicts 
  (\ref{1/3.1/2.lem}).

  If  $s$  is a smooth point of $\phi(C)$ then $\deg C\leq 2$ by (\ref{sm.at.s.lem}).

  We are left with the case when  $s\notin \phi(C)$.
Set $c':=N\cap \phi(C)$.  If $c'$ is a smooth point of $\phi(C)$ then again $\deg C\leq 2$ by (\ref{sm.at.s.lem}).
The last case is when $c'$ is a smooth point of $N $, but a singular point of $C'$.  Then $c'$ is a maximal multiplicity point of $C'$, which
 is impossible by (\ref{smooth.cont.lem}). 
 \qed

\medskip

    \begin{thm} \cite{MR1016092} \label{1/3.1/2.lem} Let $C\subset \p^2$ be a rational curve, all of whose singularities are cusps. Then   there is a point $c\in C$ such that $\mult_cC>\frac13 \deg C$. \qed
    \end{thm}

\cite{MR1942244} proves a form where $\frac13$ is  asymptotically 
   replaced by  $\frac{2}{3+\sqrt{5}}\approx \frac{1}{2.62}$.

 \begin{lem} \label{sm.at.s.lem}
Let $C\subset \p^2$ be a cuspidal, rational curve  such that
  $C\cap N =\{s\}$. If $s$ is a smooth point of $C$ then $\deg C\leq 2$.
\end{lem}

Proof.
The cases $\deg C\geq 4$  are excluded by (\ref{meet.smooth.lem}).
If  $\deg C= 3$, then $N $ lies in the pencil   $|C, 3L|$ by  (\ref{meet.smooth.lem}.2).  A cuspidal cubic has a unique flex, so this pencil is
$\lambda (y^2z-x^3)+\mu z^3=0$.
It has only 1 member that is singular at $s$, namely $( z^3=0)$. \qed

\begin{lem}\label{meet.smooth.lem}
  Let $C\subset \p^2$ be a cuspidal, rational curve of degree $d$, and
  $B\subset \p^2$ a (possibly reducible) curve. Assume that
  $c:=C\cap B$  is a  single smooth point of $C$ with tangent line $L$. Then one of the following holds.
  \begin{enumerate}
  \item  $\deg B> \deg C$,
  \item  $\deg B= \deg C$  and $B\in |C, dL|$, or
  \item  $B=eL$  for some $e<d$.
    \end{enumerate}
\end{lem}

Proof. Set  $e:=\deg B$ and let $H$ be the line class. By assumption
$de[c]$ and $e[H|_C]$ are linearly equivalent on $C$. Since $C$ is  cuspidal and rational,
$\pic(C)$ is torsion free; see, for example,  \cite[Exrc.II.6.9]{hartsh}. Thus
$d[c]\sim [H|_C]$.  That is, $L$  meets
$C$ with multiplicity $d$.   If $e<d$ then
$H^0\bigl(\p^2, \o_{\p^2}(B-C)\bigr)=0$, thus $B$ is unique, so $B=eL$.

If $e=d$ then $H^0\bigl(\p^2, \o_{\p^2}(B-C)\bigr)=1$, so
$B$ lies in the pencil  $|C, dL|$.
\qed

\begin{lem}\label{smooth.cont.lem} Let $C'\subset \p^2$ be a cuspidal, rational curve. Assume  that
  $c':=C'\cap N $ is a single smooth point on $N $ and 
   a maximal multiplicity point of $C'$. Then $C'$ is a line or a conic.
\end{lem}

Proof.  Set $r':=\deg C'$ and let 
$C'$ have a $(p',q')$-cusp at $c'$. Note that $q'>\frac13 r'$ by (\ref{1/3.1/2.lem}). 
  As we noted in (\ref{cusp.defn}),
$p'\geq(C\cdot N )= 3r'$.
 Using   (\ref{genus.say}.4) we get  that
$$
r'^2-\tfrac73 r' +\tfrac23=(3r'-1)\bigl(\tfrac13 r'-\tfrac23\bigr)\leq (r'-1)(r'-2),
$$
thus $\tfrac 23 r'\leq  \tfrac 43$. \qed

  \begin{exmp}\label{3.smooth.lines} There are 4 types of smooth, rational plane curves $C$ such that
    $C\cap N $ is a single point.
    Here $N :=(xyz=x^3+y^3)$ with parametrization
    $t\mapsto (t{:} t^2{:} 1+t^3)$.

    (\ref{3.smooth.lines}.1)  $C_1:=(x=0)$ and $(y=0)$, meeting $N $ at the node. These have  contact orders $(2,1)$ or $(1,2)$ with $N$.
    If $C$ is a 
conic with contact orders $(5,1)$, then $\sigma_{-}(C)$ is a 
line with contact orders $(2,1)$.

These conics are $(x-y^2=0)$ and $(y-x^2=0)$

 (\ref{3.smooth.lines}.2) Flex tangents.
To find these we are looking for a line $ax+by+z=0$ such that 
   $at+bt^2+(1+t^3)=(\mbox{cube})$.
The unique real  solution is  with $(1+t)^3$, giving
  $C_2:=(3x+3y+z=0)$, and   $C_2\cap N=(1{:}{-}1{:}0)$. Acting by $\rho_3$ gives 2 more.

(\ref{3.smooth.lines}.3)  5-tangent conics.
We need to solve
  $$
  at^2+bt^3+ct^4+dt(1+t^3)+et^2(1+t^3)+(1+t^3)^2=(\mbox{6th power}),
  $$
  but we want to avoid the 3 cases that give  the square of a flex tangent equation. 
  So the  unique real  solution is with $(1-t)^6$, giving
  $C_3:=\bigl(Q(x,y,z)=0\bigr)$ where
  $$
  Q(x,y):=21x^2-22xy+21y^2-6xz-6yz+z^2,
  $$
  and   $C_3\cap N=(1{:}1{:}2)$.
  See  \cite[p.670]{MR1942244} for a different derivation of this.
  Acting by $\rho_3$ gives 2 more.
  \end{exmp}

In order to complete the proof of Corollary~\ref{thm.aff.lines} it remains  to show 
that the curves $C_1, C_2, C_3$ are in different $\aut(\p^2\setminus N)$-orbits.

  \begin{say}\label{map.to.gm.say}  Let $N \neq C:=(g=0)\subset \p^2$ be an irreducible curve of degree $d$. The group of non-constant units on $\p^2\setminus(N \cup C)$
    is generated by   $(xyz-x^3-y^3)^d/g^3$  if $3\nmid d$, and by
    $(xy-x^3-y^3)^{d/3}/g$ if $3\mid d$.  This gives  a morphism
    $$
    m_C: \p^2\setminus(N \cup C)\to \p^1\setminus\{0,\infty\}
    $$
    that is canonically associated to $N $ and $C$.
We describe the
fibers  in the cases (\ref{3.smooth.lines}.1--3).

    (\ref{map.to.gm.say}.1) For $C_1$ the fibers are  $\lambda (xyz-x^3-y^3)+\mu x^3=0$; these are rational curves.

    (\ref{map.to.gm.say}.2) For $C_2$ the fibers are   $\lambda (xyz-x^3-y^3)+\mu (3x+3y+z)^3=0$, these are genus 1 curves.

    (\ref{map.to.gm.say}.3) For $C_3$ the fibers are   $\lambda (xyz-x^3-y^3)^2+\mu Q^3=0$.
    At the base point $c$
we can choose local coordinates such that $N =(u=0)$ and $Q=(v-u^5=0)$. 
    Thus the fibers have a singularity with equation  $\lambda u^2+\mu (v-u^5)^3=0$, which is
    analytically equivalent to $u^2+v^{15}=0$.  By (\ref{genus.say}.3)  the genus of this singularity is $7$, thus the geometric genus of the  general fiber is
    $\binom{5}{2}-7=3$.

    In particular, the pairs  $\bigl(\p^2\setminus N , C_i\bigr)$
    are not isomorphic to each other for different values of $i\in\{1,2,3\}$.
  \end{say}

 \begin{exmp}[Cubics with contact order  $(1,8)$]\label{3.smooth.cubics}
  For these we need to solve
  $$
  g_3(t, t^2)+g_2(t, t^2)(1+t^3)+g_1(t, t^2)(1+t^3)^2+(1+t^3)^3=t^8,
  $$
  where $g_i$ is homogeneous of degree $i$. The solutions form a pencil
  $$
  \lambda(xyz-x^3-y^3)+\mu(-xy^2-x^2z+yz^2)=0,
  \eqno{(\ref{3.smooth.cubics}.1)}
  $$
  whose general member is  smooth.

  Remarkable properties of this pencil are studied in \cite{moschetti2024pencilsplanecubicsbase}. 
 \end{exmp}

   \begin{exmp}[Cubics with contact order  $(1,7)$]\label{3.smooth.cubics.7}
  For these we need to solve
  $$
  g_3(t, t^2)+g_2(t, t^2)(1+t^3)+g_1(t, t^2)(1+t^3)^2+(1+t^3)^3=at^7+bt^8,
  $$
  where $g_i$ is homogeneous of degree $i$.
Besides the pencil   (\ref{3.smooth.cubics}.1), we also have $y(yz-x^2)$.
  The resulting linear system
  $$
  \lambda(xyz-x^3-y^3)+\mu(-xy^2-x^2z+yz^2)+\nu (-x^2y+y^2z)=0.
  $$
   gives the degree 2 map
  $$
  \p^2=S\stackrel{\pi_+}{\leftarrow} T_+
  \stackrel{p}{\rightarrow} \p^2
  $$
  in the definition  (\ref{a.obu}) of the 
Geiser involutions    for $\bigl(\p^2, (xyz=x^3+y^3)\bigr)$.
 \end{exmp}

  \begin{exmp} \label{a.l.F1.exmp.1} Let $S$ be a smooth Del~Pezzo surface and $C\subset S$ a smooth, rational curve with $(C^2)\geq 0$. Set $L:=\o_S(C)$.

    Let $N\in |-K_S|$ be a nodal curve. Consider the sequence
    $$
    0\to \omega_S(L)\to L\to L|_N\to 0.
    $$
    By the adjunction formula,  $\deg L|_N=2+(C^2)$. Thus there are
    $2+(C^2)$  sections $s_i\in H^0(N, L|_N)$ such that $(s_i=0)$ is a single smooth point  $c_i\in N$.
    Since $H^1\bigl(S, \omega_S(L)\bigr)=0$, the $s_i$  lift back to sections
    of $L$. Thus we get divisors  $C_i\in |C|$ such that $C_i\cap N=\{c_i\}$.

  \end{exmp}

  \begin{exmp} \label{a.l.F1.exmp.2} Let $N \subset \p^2$ be a nodal cubic
    and $x\in N $ a smooth point.
    Blow up $x$ to get  $S:=B_x\p^2$ and  $N\subset |-K_S|$.
    Let $E\subset S$ be the $(-1)$-curve and $F\subset S$ a fiber of the projection $\pi:S\to \p^1$ from $x$.

    For any $r\geq 1$, the general members of  $|E+rF|$ are smooth rational curves of self-intersection $2r-1$. Thus by (\ref{a.l.F1.exmp.1}) there are $2r+1$ different points
    $c_{i,r}\in N$ and  divisors  $C_{i,r}\in |E+rF|$ such that $C_{i,r}\cap N=\{c_{i,r}\}$.

    If $C_i$ is reducible, then it is of the form  $F_{i,r}+C'_{i,r}$ where
    $F_{i,r}$ is the fiber through $c_{i,r}$. Thus
    $2[c_{i,r}]\sim F|_{N}$. Subtracting this from
    $(2r+1)[c_{i,r}]\sim (E+rF)|_{N}$ we get that
    $[c_{i,r}]\sim [E|_{N}]=[E\cap N]$.
So  $E\cap N$ is a branch point of the projection of $\pi|_N$.
    This means that $x\in N $ is a flex point.
    Then  $C_{i,r}=E+rF_{i,r}$ is reducible.
    Thus we conclude the following.

    For each $r\geq 1$, the curves $C_{i,r}\in |E+rF|$ give
    \begin{enumerate}
      \item  $2r+1$ different $\aut(S\setminus N)$-orbits of affine lines on $S\setminus N$ if $x\in N $ is not a flex point, but only
      \item  $2r$ different $\aut(S\setminus N)$-orbits of affine lines  if $x\in N $ is  a flex point.
    \end{enumerate}
    \end{exmp}


\def\cprime{$'$} \def\cprime{$'$} \def\cprime{$'$} \def\cprime{$'$}
  \def\cprime{$'$} \def\dbar{\leavevmode\hbox to 0pt{\hskip.2ex
  \accent"16\hss}d} \def\cprime{$'$} \def\cprime{$'$}
  \def\polhk#1{\setbox0=\hbox{#1}{\ooalign{\hidewidth
  \lower1.5ex\hbox{`}\hidewidth\crcr\unhbox0}}} \def\cprime{$'$}
  \def\cprime{$'$} \def\cprime{$'$} \def\cprime{$'$}
  \def\polhk#1{\setbox0=\hbox{#1}{\ooalign{\hidewidth
  \lower1.5ex\hbox{`}\hidewidth\crcr\unhbox0}}} \def\cdprime{$''$}
  \def\cprime{$'$} \def\cprime{$'$} \def\cprime{$'$} \def\cprime{$'$}
\providecommand{\bysame}{\leavevmode\hbox to3em{\hrulefill}\thinspace}
\providecommand{\MR}{\relax\ifhmode\unskip\space\fi MR }
\providecommand{\MRhref}[2]{%
  \href{http://www.ams.org/mathscinet-getitem?mr=#1}{#2}
}
\providecommand{\href}[2]{#2}

  \bigskip

  Princeton University, Princeton NJ 08544-1000, \

  \email{kollar@math.princeton.edu}

\end{document}